\documentclass[12pt]{article}
\usepackage{amssymb} 
\begin{document} 

\title{Diophantine non-integrability of 
a third order recurrence with the Laurent property} 
 
 
\author{ 
A.N.W. Hone\thanks{  
Institute of Mathematics, Statistics \& 
Actuarial Science,  
University of Kent,  
Canterbury CT2~7NF, UK}  
} 

\def\underset#1#2{\mathrel{\mathop{#2}\limits_{#1}}}  
\newcommand{\haf}{{\hat{f}}}  
\newcommand{\beq}{\begin{equation}} 
\newcommand{\eeq}{\end{equation}} 
\newcommand{\bea}{\begin{eqnarray}} 
\newcommand{\eea}{\end{eqnarray}} 
\newcommand\la{{\lambda}} 
\newcommand\ka{{\kappa}} 
\newcommand\al{{\alpha}} 
\newcommand\be{{\beta}} 
\newcommand\de{{\delta}}
\newcommand\si{{\sigma}} 
\newcommand\lax{{\bf L}} 
\newcommand\mma{{\bf M}}  
\newcommand\ctop{{\mathcal{T}}}  
\newcommand\hop{{\mathcal{H}}}  
\newcommand\ep{{\epsilon}}
\newcommand\T{{\tau}}

\newcommand{\N}{{\mathbb N}}  
\newcommand{\Q}{{\mathbb Q}}  
\newcommand{\Z}{{\mathbb Z}}  
\newcommand{\C}{{\mathbb C}}
 
\newtheorem{Pro}{Proposition} 
\newtheorem{Lem}{Lemma} 
\newtheorem{The}{Theorem} 
 
\newcommand{\rd}{{\mathrm{d}}} 
\newcommand{\sgn}{{\mathrm{sgn}}} 
\def \ring {{\cal R}}

\maketitle  
 
\begin{abstract}  
We consider a one-parameter family of third order nonlinear 
recurrence relations. 
Each member of this family 
satisfies the singularity confinement test, has a 
conserved quantity, and moreover 
has the Laurent property: all of the iterates are Laurent 
polynomials in the initial data. However, 
we show that these recurrences 
are not Diophantine integrable according 
to the definition proposed by Halburd  
(2005 {\it J. Phys. A: Math. Gen.} {\bf 38} 
L1).  
Explicit bounds 
on the asymptotic growth of the heights of 
iterates are obtained for a special choice of 
initial data. 
As a by-product of our analysis, 
infinitely many solutions are found for a certain family  
of Diophantine equations, studied by Mordell, 
that includes Markoff's equation.   
\end{abstract}

 

For some time there has been considerable interest in 
maps or discrete equations which are integrable. 
Various different criteria have been proposed as  
tests for integrability in the discrete setting. One of 
the earliest proposals was the singularity confinement test 
of Grammaticos, Ramani and 
Papageorgiou \cite{grp}, which has proved to be an extremely useful tool 
for isolating discrete Painlev\'e equations (see e.g. \cite{rgs}). 
However, after Hietarinta and Viallet's discovery of 
some  non-integrable equations with the 
singularity confinement property, they were led to introduce 
the zero algebraic entropy condition for 
integrability of rational maps \cite{hv}, namely that 
the degree $d_n$ of the $n$th iterate of a map (as a 
rational function of the initial data)  
should satisfy $\lim_{n\to\infty}(\log d_n)/n=0$. 
The phenomenon of weak 
degree growth for integrable maps had been 
studied earlier by Veselov \cite{ves, ves2}, 
and 
the notion
of algebraic entropy proposed in \cite{hv} is connected
with other measures of growth such as Arnold complexity \cite{aabhm}.

Since the work of Okamoto \cite{oka} it has been known that 
the B\"acklund transformations of differential Painlev\'e equations 
can be classified in terms of 
affine Weyl groups. Noumi and Yamada have shown that 
the symmetries of these ODEs are simultaneously compatible with 
associated discrete Painlev\'e equations \cite{noumi}, 
and this viewpoint has been explained geometrically 
by means of Cremona group actions on rational surfaces \cite{sakai}.   
Discrete integrable systems are also related to the notions of 
discrete geometry \cite{digp} and discrete analytic functions \cite{bms}. 


In the past few years, several new techniques have been developed 
for isolating 
integrable discrete equations: analytical criteria on the one hand, and 
arithmetical ones on the other.  
The Painlev\'e property has been used very effectively to 
isolate integrable differential equations \cite{conte},  
and this led Ablowitz, Halburd and Herbst to 
extend this property to difference equations  
using Nevanlinna theory \cite{ahh}, by considering the asymptotic 
growth of meromorphic solutions at infinity. 
Roberts and Vivaldi \cite{rv} 
have studied the distribution of orbit lengths in rational 
maps reduced to finite 
fields $\mathbb{F}_p$ for different primes $p$, in order to 
identify integrable cases of such maps.   
Most recently \cite{halburd}, Halburd used Vojta's dictionary 
between Nevanlinna theory and Diophantine approximation 
to translate the concepts of \cite{ahh} 
into a Diophantine integrability 
criterion for discrete equations. 

There is by now a considerable literature on 
discrete bilinear equations, including the 
bilinear forms of discrete  
Painlev\'e  equations \cite{rgs}, and 
bilinear partial difference equations such as 
the discrete Hirota equation, which can be written  
(with two arbitrary parameters $\al ,\be$)  
in the form 
\beq 
\tau_{\ell+1,m,n}\tau_{\ell-1,m,n}
=\al \,\tau_{\ell,m+1,n}\tau_{\ell,m-1,n}
+
\be \,\tau_{\ell,m,n+1}\tau_{\ell,m,n-1}.
\label{dhir} 
\eeq 
The partial difference equation (\ref{dhir}) has   
continuum limits to all of the bilinear 
equations in the KP hierarchy of PDEs \cite{zabrodin},   
and also appears as an equation for 
transfer matrices in quantum integrable models \cite{quant}. 
However, despite the huge amount of theory that has been 
developed for bilinear discrete equations, there is another 
property of such equations which 
appears to have been overlooked by the integrable systems community, namely 
the following: for initial data 
specified on a suitable
subset of $\mathbb{Z}^3$, 
all of the iterates of (\ref{dhir})  
are Laurent polynomials in the initial data 
with coefficients in $\Z [\al ,\be ]$. 
This \textit{Laurent property}  
was originally known only to a few algebraic combinatorialists 
\cite{propp1}, and for the equation (\ref{dhir}) 
it was first proved by Fomin and Zelevinsky within the 
framework of the theory of cluster algebras \cite{fz}, 
while more detailed properties of the associated Laurent 
polynomials have been shown by Speyer \cite{speyer}. 

One of the simplest manifestations of this Laurent  phenomenon 
was found by Michael Somos, who considered bilinear 
recurrences of the form 
\beq  \label{somos}  
\tau_{n+k}\tau_n=\sum_{j=1}^{[k/2]}\tau_{n+k-j}\tau_{n+j}, \qquad k\geq 4  
\eeq   
taking the initial values $\tau_0=\tau_1=\ldots =\tau_{k-1}=1$. 
Clearly for the recurrences (\ref{somos}) each new iterate 
$\tau_{n+k}$ is a rational function of the initial data, so 
with this particular choice one would expect $\tau_n$ to be rational 
numbers, but it was observed numerically that for the Somos-4 recurrence 
\beq
\tau_{n+4}
\tau_{n}
=\al \,
\tau_{n+3}
\tau_{n+1}
+\beta
\,
(\tau_{n+2})^2 
\label{bil}
\eeq
with parameters $\alpha=\beta=1$ and starting with four ones, an 
integer sequence 
\beq \label{s4seq}
1,1,1,1,2,3,7,23,59,314,1529, 8209, 83313,\ldots
\eeq
results. Similar empirical observations showed that Somos-5, -6 and -7 
also yield integer sequences, 
while the Somos-$k$ recurrences for $k\geq 8$ do not. 
The first inductive proof of the integrality 
of the sequence (\ref{s4seq}) appeared in \cite{gale}, but it was 
realized that the deeper reason behind this lay in the fact that the   
recurrence (\ref{bil}) has the Laurent property, that is 
$\tau_n\in \mathbb{Z}[\al , \be , \tau_0^{\pm 1},
\tau_1^{\pm 1},  \tau_2^{\pm 1},  \tau_3^{\pm 1}]$ for all $n\in \mathbb{Z}$.
In other words, the iterates of (\ref{bil}) are Laurent 
polynomials in the four initial data whose coefficients are in 
$\Z [\al ,\be ]$, so if $\al =\be = \tau_0 =\tau_1=\tau_2=\tau_3=1$ 
then the sequence must consist entirely of integers. 

Fomin and Zelevinsky found a suitable setting within which 
to prove
the Laurent property for a variety 
of discrete equations \cite{fz}, including  certain recurrences of the form 
\beq \label{frec} 
\tau_{n+k}\tau_n=F(\tau_{n+1},\ldots,\tau_{n+k-1}) 
\eeq 
for particular polynomials $F$ (mostly, but not all, 
quadratic forms in their arguments), as well as integrable two- and 
three-dimensional bilinear recurrences like (\ref{dhir}).   
More recently, the connection between the iterates of Somos-4 
recurrences and sequences of points on 
elliptic curves was explained in the PhD thesis of Swart 
\cite{swart}, while independently the author found 
the explicit solution of the initial value problem 
for both Somos-4 \cite{honeblms} and Somos-5 \cite{honetams} in terms 
of elliptic sigma functions; see also \cite{vdp} for a connection 
with continued fractions. An essential observation in 
\cite{honeblms} and \cite{honetams} was that each of these 
bilinear recurrences (fourth and  
fifth order respectively) could be understood in terms of a suitable 
integrable mapping of the plane, with each one corresponding 
to a particular degenerate case of the family of maps studied by 
Quispel, Roberts and Thompson \cite{qrt}. Somos sequences are of 
considerable interest to number theorists due to 
the way that new prime factors appear therein \cite{ems, recs}.

A natural question raised by all of the above is the following: if  
discrete equations with the Laurent property 
like (\ref{dhir}) and (\ref{bil}) are 
integrable,   
then are all discrete 
equations with the Laurent property integrable? The converse 
is clearly not possible: most integrable rational maps 
do not have iterates that are Laurent polynomials in the 
initial data. Thus it is sensible to restrict ourselves 
to a suitable sub-class of maps, namely those defined by 
recurrences of the form (\ref{frec}). The purpose of this note 
is to show that the answer to this question is a negative one, 
by taking an example not considered in \cite{fz}, 
namely the one-parameter family 
of third order recurrences given 
by 
\beq 
\tau_{n+3}\tau_n=\tau_{n+2}^2+\tau_{n+1}^2 +J.  
\label{jrec} 
\eeq 
In the case $J=0$,  it was noticed by 
Dana Scott \cite{gale} 
that with $\tau_0=\tau_1=\tau_2=1$ an integer sequence 
results: 
\beq \label{markoff} 
1,1,1,2,5,29,433,37666,48928105,5528778008357,\ldots 
\eeq 
In fact, (\ref{jrec}) is related to some old 
problems in number theory: the numbers (\ref{markoff}) 
are a subsequence of the Markoff numbers that arise in the 
theory of indefinite quadratic forms \cite{markoff}, 
while the case of arbitrary $J$ is related to 
a Diophantine equation considered by Mordell \cite{mordell}. 
Below we give two different proofs that the recurrence (\ref{jrec}) 
has the Laurent property; 
the first proof also shows 
that it passes the singularity confinement test, while the second 
one gives a stronger result based on the existence of a conserved quantity.     
However, when Halburd's Diophantine integrability test is applied to 
a particular class of sequences generated by (\ref{jrec}), 
it is seen that the growth 
of the heights of iterates is double exponential (which can already 
be guessed from the rapid growth of the terms in (\ref{markoff})). For the  
the particular sequences being considered,  quite sharp bounds on 
this double exponential growth are obtained.  
Hence 
we have an interesting example of a non-integrable recurrence 
with the Laurent 
property. 
 

The recurrence (\ref{jrec}) defines a 
rational  map of three-dimensional affine space, i.e.  
\beq 
\label{3d} 
(x_n,y_n,z_n) \mapsto 
(y_n , z_n , (y_n^2+z_n^2+J)/x_n ) 
= 
(x_{n+1},y_{n+1},z_{n+1}) , 
\eeq 
upon setting
$(\tau_n,\tau_{n+1},\tau_{n+2})=(x_n,y_n,z_n)$. 
Before considering the properties of this map in more detail, 
we present the first proof that it has the Laurent property, since 
this will simultaneously show that it satisfies singularity confinement 
as well. 
Letting $\tau_0=a$, $\tau_1=b$, $\tau_2=c$ and putting the recurrence into 
MAPLE, we see that the next two iterates are  
$$ 
\tau_3=\frac{b^2+c^2+J}{a}, \quad 
\tau_4=\frac{b^4+a^2c^2+2b^2c^2+c^4+(a^2+2b^2+2c^2)J+J^2}{a^2b}, 
$$  
and $\tau_5$ has $a^4b^2c$ in the denominator and its numerator is 
a polynomial in $\Z [a,b,c,J]$ with 29 terms. The first 
miracle occurs when $n=3$, because upon dividing 
the right hand side of (\ref{jrec}) by $\tau_3$ the numerator $b^2+c^2+J$ 
cancels, so that 
$\tau_6\in \ring =\Z [a^{\pm 1},b^{\pm 1},c^{\pm 1},J]$ 
with $a^7b^4c^2$ in the denominator and a polynomial with 
104 terms in the numerator.  

Now to prove the Laurent property, 
we can use the fact that modulo 
monomial factors $a^jb^kc^\ell$  the ring $\ring$ of Laurent polynomials 
is a unique factorization domain (see \cite{fz} - 
the case $J=0$ is covered by Theorem 1.10). Thus we can consider   
divisibility of the $\tau_n$ modulo such monomials. As the  
inductive hypothesis, assume that $\tau_j\in\ring$ for $0\leq j\leq n+5$, 
that any three adjacent terms $\tau_j,\tau_{j+1},\tau_{j+2}$ 
are pairwise coprime (modulo monomials) in this range, and 
also that $\tau_{j\pm 1}^2+J$ is coprime to $\tau_j$. To 
prove that the next iterate $\tau_{n+6}$ is in $\ring$, 
consider 
the three adjacent Laurent polynomials $x=\tau_{n+1}$, 
$y=\tau_{n+2}$, $z=\tau_{n+3}$. From successive applications 
of (\ref{jrec}) we have 
\beq   
\label{smodz} 
x^2+y^2+J=\tau_{n+3}\tau_n\equiv 0 \quad \bmod z, 
\eeq 
as well as  
$ 
\tau_{n+4}\equiv (y^2+J)/x
$ and  
$ 
\tau_{n+4}\equiv (y^4+J(x^2+2y^2)+J^2)/(x^2y) 
$ 
(both $\bmod z$), 
which together imply that 
\beq 
\label{finalt} 
\tau_{n+6}\tau_{n+3}\equiv \frac{(y^2+J)(x^2+y^2+J)(y^4+J(x^2+2y^2)+J^2)} 
{x^4y^2}\equiv 0 \, \bmod z 
\eeq 
by (\ref{smodz}), and hence $\tau_{n+6}\in\ring$. 
Next suppose that an irreducible factor 
$p\in\ring$ is such that $p|\tau_{n+6}$ and 
$p|\tau_{n+5}$; then $p|(\T_{n+6}\T_{n+3}-\tau_{n+5}^6)=\T_{n+4}^2+J$, 
which contradicts the coprimality of $\T_{n+5}$ and 
$\T_{n+4}^2+J$, so $\T_{n+6}$ and $\T_{n+5}$ must be coprime. 
The same argument shows that 
$\T_{n+6}$ and $\T_{n+4}$ must also be coprime, so that the 
three adjacent terms $\T_{n+4}$, $\T_{n+5}$ and 
$\T_{n+6}$ are pairwise coprime. Similarly it follows that 
$\T_{n+5}^2+J$ is coprime to $\T_{n+6}$, and also 
$\T_{n+6}^2+J$ is coprime to $\T_{n+5}$, which completes 
the inductive step. This proves that $\T_n\in\ring$ for 
all positive values of the index $n$, and the same conclusion holds for 
negative indices by the reversibility of the recurrence. 

To see why singularity confinement also follows from the  preceding 
argument, observe that a singularity can only occur when one of 
the iterates vanishes. Thus, without loss of generality 
suppose that $z=\T_{n+3}$ 
vanishes, and we take $\T_n=w$ and $y^2=-x^2-J+\ep w$  
so that $z=\ep\to 0$. The above expressions can quickly be evaluated in this 
limit by noting that all terms congruent to $z$ are $O(\ep )$, and 
also 
$y=y_0+O(\ep )$ with $y_0=\sqrt{-x^2-J}$. Thus $\T_{n+4}=-x+O(\ep )$, 
$\T_{n+5}=(x^2+J)/y_0+O(\ep )$, and clearly from (\ref{finalt}) 
we have $\T_{n+6}\,\ep =O(\ep )$ so $\T_{n+6}=O(1)$ and the singularity is 
confined. This argument for confinement will work for many other 
recurrences with the Laurent property, such as those in \cite{fz}. For 
the recurrence (\ref{jrec}), a  
stronger version of the Laurent property actually holds, as a 
consequence of the following

\noindent \textbf{Proposition 1} \textit{ 
The quantity 
$N_n=(\T_{n+2}^2+\T_{n+1}^2+\T_{n}^2+J)/(\T_{n+2}\T_{n+1}\T_{n})$ 
is an invariant for the recurrence (\ref{jrec}), so that 
$N_n=N_0=N$ with 
\beq 
N=\frac{a^2+b^2+c^2+J}{abc}. 
\label{ne} 
\eeq    
Moreover, for the same  initial data $\T_0=a, \T_1=b, \T_2=c$, 
the iterates of (\ref{jrec}) are the same as those of the recurrence 
\beq 
\label{nrec} 
\T_{n+3}=N\T_{n+2}\T_{n+1}-\T_{n}, 
\eeq 
and all $\T_n$ lie in the ring  
$\Z [a,b,c,N]\subset \ring=\Z [a^{\pm 1},b^{\pm 1},c^{\pm 1},J]$. 
} 
  
The proof is quite straightforward. The fact that $N_n$ is a 
conserved quantity for (\ref{jrec}) 
follows from the  direct calculation that 
$N_{n+1}-N_n=0$, so clearly the iterates of the map (\ref{3d}) 
lie on the cubic surface 
\beq 
\label{surf} 
x^2+y^2+z^2-Nxyz+J=0, 
\eeq 
which is non-singular for $J\neq 0,-4/N^2$.    
This surface admits some obvious symmetries: it is unchanged  
by permutation of $(x,y,z)$, and also for fixed $y$ and $z$ it is 
a quadratic equation in $x$ and hence is preserved by the 
involution $\mathcal{I}:\, (x,y,z)\mapsto (x^\dagger ,y,z)$ 
where  
\beq 
x^\dagger =(y^2+z^2+J)/x=Nyz-x 
\label{dag} 
\eeq   
is the other root of the quadratic. 
Taking the cyclic permutation 
$\mathcal{P}:\,(x,y,z)\mapsto (y,z,x)$, 
then performing $\mathcal{I}$ 
followed by 
$\mathcal{P}$ yields the transformation 
$(x,y,z)\mapsto (y,z,x^\dagger )$. 
If the first equality for $x^\dagger$ in (\ref{dag})
is used in this transformation, then the  rational map (\ref{3d}) 
corresponding to the recurrence (\ref{jrec}) arises. 
Similarly, the 
alternative recurrence (\ref{nrec}) comes  from  
the second equality  in (\ref{dag}), and one could start with 
this and observe that it has the conserved quantity  
$J_n =N\T_{n+2}\T_{n+1}\T_n-\T_{n+2}^2-\T_{n+1}^2-\T_{n}^2$. 
It is also obvious 
that the iterates of (\ref{nrec}) lie in  
$\Z [a,b,c,N]$, and with $N$ given by the formula 
(\ref{ne}) this is a subring of $\ring$.    

Here we should point out that the aforementioned 
symmetries 
were used extensively by Mordell \cite{mordell}  when he  
considered the Diophantine problem of finding integer triples  
$(x,y,z)$ satisfying (\ref{surf}) for fixed integers $N$ and $J$. 
The case $N=-2$, $J=-n\in\Z$ was treated in great detail, since this is 
related to the problem of finding which integers  can be written 
as the sum of four cubes. However, he remarked that 
``\textit{I know of no test except trial for the existence of solutions 
for given $n$.}'' The particular 
case $N=3$, $J=0$ is another distinguished 
instance of (\ref{surf}) 
known as Markoff's equation, which is related to the spectrum 
of indefinite binary quadratic forms \cite{markoff}. Although a great 
deal more is known about this case, there are still difficult open 
problems associated with it \cite{baragar, button, perrine}.  
The
case $J=0$ has also acquired a combinatorial interpretation 
very recently 
\cite{itsara, propp2}.
The above Proposition implies that, 
given any initial solution triple of integers  
$(a,b,c)$, the recurrence (\ref{jrec}) (or (\ref{nrec}), which 
is equivalent) 
will generate infinitely many such triples provided that it  
leads to neither a zero nor a periodic orbit. In fact, by the 
above discussion of singularity confinement, a single zero (e.g. a triple 
$(x,y,0)$) is not a problem in the sense that one can introduce 
a small parameter $\ep$ and analytically continue through the apparent 
singularity. However, a solution like $(x,0,0)$ can also 
occur if $-J$ is a perfect square: see the Lemma and 
the Theorem on p.506 of 
\cite{mordell}; this is a more serious obstacle. 
On the other hand, given an initial solution 
$(x_0,y_0,z_0)$ of (\ref{surf}),  
by that Lemma one can still generate infinitely many solutions 
$(x_n,x_{n+1},z_0)$ with the same 
value of $z$ by setting $x_1=y_0$ and iterating the  
\textit{integrable} recurrence 
\beq \label{integ} 
x_{n+1}x_{n-1}=x_n^2+z_0^2+J  
\eeq 
(with $N^2z_0^2\geq  4$ for an infinity of integer solutions) 
or, equivalently, the linear recurrence 
\beq 
\label{linear} 
x_{n+1}-Nz_0x_n+x_{n-1}=0. 
\eeq


Halburd's Diophantine integrability test \cite{halburd} 
applies to rational maps 
with iterates and parameters taking values in $\Q$ or a number field. 
For $x\in\Q$ with $x=p/q$ as a fraction in lowest terms, the height 
is $H(x)=\max \{ |p|,|q| \}$, while the logarithmic height is 
$h(x)=\log H(x)$. A rational map of this kind 
is defined to be \textit{Diophantine integrable} 
if its iterates $x_n$ are  such  
that  
$h(x_n)$  
grows no faster than a polynomial in $n$. To apply the test to the recurrence 
(\ref{jrec}) we shall restrict ourselves to the case when all 
of the iterates $\tau_n$ are integers, so that 
$h(\tau_n)=\log |\T_n |$. More generally we can consider 
real or complex values of $\T_n$, and  
set $\Lambda_n=\log |\T_n |$. Now suppose that 
$\Lambda_n\sim C\la^n\to\infty$ as $n\to\infty$ 
for real $\la >1$ and some $C>0$. Then clearly $\T_n^{-1}\to 0$ 
and also $\Lambda_{n+1}-\Lambda_{n+2}\sim C\la^{n+1}(1-\la )\to -\infty$, 
whence $\T_{n+1}/\T_{n+2}\to 0$ as $n\to\infty$. So taking the  
logarithm of (\ref{jrec}) and substituting in these asymptotics 
gives 
\beq 
\label{asy} 
\Lambda_{n+3}-2\Lambda_{n+2}+\Lambda_n=
\log\left|1+\frac{\T_{n+1}^2}{\T_{n+2}^2}+\frac{J}{\T_{n+2}^2}\right| 
\to 0. 
\eeq 
The characteristic polynomial for the linear recurrence on the 
left hand side of (\ref{asy}) gives $\la^3-2\la^2+1=0$, with the 
largest root being the golden mean, $\la=(1+\sqrt{5})/2>1$, 
which is 
consistent with the original assumptions on the asymptotics of 
$\Lambda_n$. Thus if integer sequences of this kind exist, then 
their logarithmic heights grow exponentially with $n$ with 
$\lim_{n\to\infty} (\log h(\tau_n))/n = \log ((1+\sqrt{5})/2)\approx 
0.4812$. We proceed to show existence. 
          
\vspace{0.05in} 
\noindent \textbf{Proposition 2} 
\textit{ 
With initial 
data 
$(\T_0,\T_1,\T_2)=(1,1,1)$  
the recurrence (\ref{jrec}) gives a sequence of monic 
polynomials 
$\T_n=p_n(N)\in\Z [N]$ with $N=J+3$, and  
$\deg p_n=f_n-1$ for $n>0$ where $f_n$ is 
a Fibonacci number. For real $N>2$, 
$p_{n+1}(N)>p_{n}(N)>1$  holds for $n\geq 3$. 
Furthermore, for 
$N>2$ and all $n\geq 4$
these polynomials satisfy
\beq \label{bounds}
(N-1)^{f_n-1}<p_n(N)<N^{f_n-1}.
\eeq 
} 

That these initial data give polynomials in $N$ follows from the 
earlier Proposition; equation (\ref{ne}) gives $N=J+3$,
and we have $p_3=N-1$, 
$p_4=N^2-N-1$, $p_5=N^4-2N^3+N-1$, $p_6=N^7-3N^6+N^5+3N^4-2N^3+1$ etc. 
From (\ref{nrec}), $d_n=\deg p_n$ satisfies 
$d_{n+3}=d_{n+2}+d_{n+1}+1$ which implies that $d_n=f_n-1$,  
with $f_n$ being a Fibonacci number ($f_0=0$, $f_1=1$), and 
the $p_n$ are clearly monic.  
The values $p_n(-1)$, $p_n(0)$ and $p_n(1)$ cycle with periods 4,6 and 
12 respectively, while 
$p_n(2)=1$ for all $n$, but for $N>2$ we have $J>-1$ 
so by induction $p_{n+1}^2+J>0$ and (\ref{jrec}) gives 
$p_{n+3}p_n>p_{n+2}^2$ so  
$p_{n+3}/p_{n+2}>(p_{n+2}/p_{n+1})(p_{n+1}/p_{n})>1$ as required. 
The upper bound in (\ref{bounds}) 
is easy to see by induction from (\ref{nrec}), 
since $p_{n+3}<Np_{n+2}p_{n+1}$. The lower bound follows from 
another induction, by writing  
$p_{n+1}=p_{n}\rho_n$, and noting 
that for $N>2$ the ratios satisfy  
$\rho_{n}>(N-1)^{f_{n-1}}$ 
when $n\geq 3$, as a consequence of the inequality 
$\rho_{n+2}>\rho_{n+1}\rho_{n}$. 
Our main result is an immediate consequence 
of the preceding facts. 

\vspace{0.05in}
\noindent \textbf{Theorem}
\textit{For each integer $N\geq 3$, every triple   
$(p_{n}(N),p_{n+1}(N),p_{n+2}(N))$ of adjacent 
iterates in the sequence for $n\geq 0$ constitutes a distinct 
solution triple for the Diophantine equation 
(\ref{surf}) with $J=N-3$. 
Moreover, for all such $N\in\Z$ 
the logarithmic height $h(p_n(N))$ of these iterates  
grows exponentially, and for all real 
$N>2$ we have 
$\lim_{n\to\infty}(\log\log p_n(N))/n =   
\log \la \approx 0.4812$ where $\la$ is the golden mean. 
} 

The integrality of the terms $p_n(N)$ for $N\in\Z$ is obvious. 
The fact that these are distinct solution triples follows 
from the monotonicity of the sequence for $N>2$, while the 
statements about asymptotic growth follow from taking logarithms 
of each side of (\ref{bounds}). 
Similar estimates hold for $|p_n|$ when 
$N\leq  -2$, but then $p_{3n}(N)<0$. It is quite easy to see that $\log\la$ 
is also the value of the algebraic entropy for the recurrence.   

It is hoped that the connections between Diophantine integrability, 
singularity confinement and the Laurent phenomenon will become 
clearer in the future.        
\vspace{.05in} 

\small 
\noindent {\bf Acknowledgements.} 
The author is grateful for the support of an EPSRC Springboard 
Fellowship, and thanks Jack Button, Nalini Joshi, Yuri Suris and 
Sasha Veselov for useful 
discussions. 


\end{document}